\documentclass[letterpaper, reqno,12pt]{amsart}
\usepackage[utf8x]{inputenc}
\usepackage{amsmath, amssymb, amsfonts}
\usepackage[foot]{amsaddr}
\usepackage{amsthm}

\usepackage{enumerate}
\usepackage[margin=2.5cm]{geometry}

\title{Asymptotic Behavior of the Expectation Value of Permanent Products}

\author{Paul Federbush$^*$}
\address{$^*$Department of Mathematics, University of Michigan, Ann Arbor, MI 48109-1043, \textit{email:} pfed@umich.edu}

\date{July 21, 2014}

\begin{document}

\maketitle

\begin{abstract}
    We would desire to have done the calculations of this paper in the measure on $n \times n$ matrices that weights uniformly all $0-1$ matrices with row and column sums equal to $r$, other matrices given weight zero. Instead we work with all matrices that are the sum of $r$ independent uniformly weighted permutation matrices, with the hope that the computations we perform give the same result in this measure. We derive the result for limiting expectations
    $$\lim_{n \to \infty} \frac{1}{n} \ln( E( \text{perm}_m (A) \: \text{perm}_{m'} (A) ) ) = \lim_{n \to \infty} \frac{1}{n} \ln( E( \text{perm}_m (A) ) )+ \lim_{n \to \infty} \frac{1}{n} \ln( E( \text{perm}_{m'} (A) ) )$$
    Here $r$ is fixed and $m$ and $m'$ are taken as each proportional to $n$.
\end{abstract}
\vfill
So we work with a measure on matrices giving uniform weight to a sum of $r$ independent permutation matrices. \cite{1}. We have as a goal the formula
\begin{equation}
    \lim_{n \to \infty} \frac{1}{n} \ln( E( \text{perm}_m (A) \: \text{perm}_{m'} (A) ) ) = \lim_{n \to \infty} \frac{1}{n} \ln( E( \text{perm}_m (A) ) )+ \lim_{n \to \infty} \frac{1}{n} \ln( E( \text{perm}_{m'} (A) ) )
\end{equation}
where $r$ is fixed, and $m$ and $m'$ are each proportional to $n$. The result may seem more reasonable in light of Friedland's Asymptotic Matching Conjecture, now proven. \cite{2, 3, 4}. Many readers will feel the computation of this paper can in a straightforward manner be made into a rigorous proof. But we do make an assumption later, the ``statistical assumption'', that though ``obvious'', should probably have a formal proof. We plan to present such a proof in a future paper.

We start with the computation of the expectation of perm$_m( A)$ given in \cite{5}, which serves as a model for us. Their equation (4.3):
\begin{equation} \label{eq:2}
    E( \text{perm}_m A ) = \frac{1}{(n!)^r} \binom{n}{m}^2 m! \sum_{m_1, \ldots, m_r \in \mathbb{Z_+}, m_1 + \ldots m_r = m} \frac{m! (n - m_1)! \ldots (n - m_r)!}{m_1! \ldots m_r!}
\end{equation}
We let $F(z)$ be given by
\begin{equation} \label{eq:3}
    F(z) = z \ln z - z
\end{equation}
We define the set $S$ to be the region
\begin{equation}
    \begin{split}
        & m_i \geq 0 \quad i = 1, \cdots, r-1 \\
        & m_1 + \cdots + m_{r - 1}\leq m
    \end{split}
\end{equation}
and set
\begin{equation} \label{eq:5}
    m_r = m - m_i - \cdots - m_{r-1}
\end{equation}
We define
\begin{equation}
    R = -r F(n) + 2 F(n) - 2 F(n - m) - \sum_i F(m_i) + \sum_i F(n - m_i)
\end{equation}
We now have as an approximation to \eqref{eq:2} the formula
\begin{equation} \label{eq:7}
    \int_S dm_1 \cdots dm_{r - 1} e^R
\end{equation}
The maximum of $R$ inside the region of integration is achieved when all $m_i$ are equal, so $m_i = \frac{m}{r}$. Let $R_0$ be this maximum value of R. Then
\begin{equation}
    \lim_{n \to \infty} \frac{1}{n} \ln( E( \text{perm}_m ) ) = \frac{R_0}{n}
\end{equation}
With
\begin{equation}
    m = pn
\end{equation}
and so $m_i = \dfrac{pn}{r}$, $\dfrac{R_0}{n}$ becomes
\begin{equation}
    -p \ln p + (2p - r) \ln r + 2 (p - 1) \ln (1 - p) + (r - p) \ln (r- p)
\end{equation}
We will play the same game with the expectation $E( \text{perm}_m \: \text{perm}_{m'} )$.

We introduce some notation. We view the sum in \eqref{eq:2} as over ``$m-terms$.'' An \underline{m-term} is an (unordered) set of m \underline{elements}, each element having a \underline{location} and a \underline{color}. The color is an integer in the set $1, 2, \cdots, r$, and the location is an assignment of a row and a column in a universal $n \times n$ matrix. There is the restriction that two distinct elements in an m-term cannot have the same row or the same column in their
locations. In \eqref{eq:2} the expression
\begin{equation}
    \binom{n}{m}^2 m!
\end{equation}
sums over the locations, and the expression
\begin{equation}
    \sum_{m_i} \frac{m!}{m_1! \cdots m_r!}
\end{equation}
sums over the assignments of color.

In computing $E( \text{perm}_m \: \text{perm}_{m'} )$ we will be summing over pairs consisting of an m-term and an m$'$-term. There is the additional restriction that if an element of the m-term is of the same color as an element of the m$'$-term, and if these two elements have the same row or column in their locations, then they have the same location (i.e. they have the same row and same column). Dealing with this restriction is the source of all the difficulties.

We construct an exact expression for the expectation of a product of permanents, a generalization of \eqref{eq:2}
\begin{equation} \label{eq:13}
    E( \text{perm}_m \: \text{perm}_{m'} ) = M \circ A \circ E \circ B \circ C \circ D \circ T
\end{equation}
Each module in the right side of \eqref{eq:13}, except the last, includes sums. The sums are to be performed from right to left, the same order as in integration.

The \underline{module $M$} sums over the m-terms in perm$_m$
\begin{equation} \label{eq:14}
    M = \sum_{m_i} \frac{1}{(n!)^r} \binom{n}{m}^2 m! \frac{m!}{m_1! \cdots m_r!}
\end{equation}
$m_i$ is the number of elements in a given m-term of color $i$. Equation \eqref{eq:5} holds. We have included the normalization factor $\dfrac{1}{(n!)^r}$ in $M$. Specification of ranges of sums is postponed to later. Note that we thus sum over the $ {m'}$-terms in perm$_{m'}$ in a ``background'' of a fixed m-term, later summed over by $M$. $A$, $E$, $B$, $C$, and $D$ sum over these m$'$-terms.

With the fixed background m-term (laid down by $M$) given, we divide the elements of any m$'$-term into five classes.
\begin{description}
    \item[class 1] An element is in class 1 if the row and column of its location are disjoint from all the rows and columns associated to elements in the fixed m-term. class 1 is summed over by $A$.
    \item[class 2] An element is in class 2 if it agrees with some element in the fixed m-term in both color and location. class 2 is summed over by $E$.
    \item[class 3] An element is in class 3 if the row of its location agrees with the row assignment of some element in the fixed m-terms, but its column assignment is disjoint from that of all elements in the fixed m-term. class 3 is summed over by $B$.
    \item[class 4] An element is in class 4 if the column of its location agrees with the column assignment of some element in the fixed m-term, but its row assignment is disjoint from that of all elements in the fixed m-term. class 4 is summed over by $C$.
    \item[class 5] An element is in class 5 if its row assignment agrees with the row assignment of some element in the fixed m-term, and likewise its column assignment agrees with the column assignment of some element in the fixed m-term, but it does not agree in both color and location with any element in the m-term. class 5 is summed over by $D$.
\end{description}

For the m$'$-term being summed up let $a$ be the number of its elements in class 1, $e$ the number in class 2, $b$ the number in class 3, $c$ the number in class 4, and $d$ the number in class 5. The classes are disjoint and cover all cases, so one has
\begin{equation} \label{eq:15}
    a + e + b + c + d = m'
\end{equation}
For each letter adding a single subscript specifies the number of a given color. So $b_i$ is the number of elements in class 3 of color $i$. One has
\begin{equation} \label{eq:16}
    a = \sum_i a_i
\end{equation}
\begin{equation} \label{eq:17}
    e = \sum_i e_i
\end{equation}
\begin{equation} \label{eq:18}
    b = \sum_i b_i
\end{equation}
\begin{equation} \label{eq:19}
    c = \sum_i c_i
\end{equation}
\begin{equation} \label{eq:20}
    d = \sum_i d_i
\end{equation}

The \underline{module $T$} sums over all possibilities for the $r$ permutation matrices not determined by the selection of the given m-term and m$'$-term.
\begin{equation} \label{eq:21}
    T = \prod_i (n - m_i - a_i - b_i - c_i - d_i)!
\end{equation}
This is the generalization of the
\begin{equation}
    \prod_i (n - m_i)!
\end{equation}
factor in \eqref{eq:2}.

The \underline{module $A$} sums over all possibilities for the number of elements in the m$'$-term of class 1.
\begin{equation} \label{eq:23}
    A = \sum_{a_i} \binom{n - m}{a}^2 a! \left( \frac{a!}{a_1! \cdots a_r!} \right)
\end{equation}

\underline{Module $E$} sums over all possibilities for the number of elements in the m$'$-term of class 2.
\begin{equation} \label{eq:24}
    E = \sum_{e_i} \prod_i \binom{m_i}{e_i}
\end{equation}

We now consider \underline{module $B$} which sums over the number of elements in the m$'$-term of class 3. We have non-negative integers $b_{ik}$ with $i, k \in \{1, 2, \cdots, r\}$ and $i \neq k$.
\begin{equation}
    b_i = \sum_{k \neq i} b_{ik} \label{eq:25}
\end{equation}
We also define
\begin{equation}
    \widetilde{b_i} = \sum_{k \neq i} b_{ki} \label{eq:26}
\end{equation}
$b_{ik}$ is the number of elements in the m$'$-term that are of class 3 and color $i$ and which share a row with an element of the m-term of color $k$.
\begin{equation} \label{eq:27}
    B = \sum_{b_{ik}} \left[ \binom{n - m - a}{b} \frac{b!}{b_1! \cdots b_r!} \prod_i \frac{b_i!}{\prod_{k \neq i} b_{ik}!} \right] \cdot \left[ \prod_i \binom{m_i - e_i}{\widetilde{b_i}} \frac{\widetilde{b_i}}{\prod_{k \neq i} b_{ki}!} \right] \cdot \left[ \prod_i \prod_{k \neq i} b_{ik}! \right]
\end{equation}
The first bracket selects the columns to be assigned to the elements of each set $b_{ik}$, the second bracket selects the rows to be assigned to the elements of each set $b_{ik}$, and the third bracket counts the ways of grouping the rows and columns.

\underline{Module $C$} is constructed as module $B$ with rows and columns interchanged. $c_i$ and $c_{ik}$ replace $b_i$ and $b_{ik}$.
\begin{equation}
    c_i = \sum_{k \neq i} c_{ik} \label{eq:28}
\end{equation}
\begin{equation}
    \widetilde{c_i} = \sum_{k \neq i} c_{ki} \label{eq:29}
\end{equation}
\begin{equation} \label{eq:30}
    C = \sum_{c_{ik}} \left[ \binom{n - m - a}{c} \frac{c!}{c_1! \cdots c_r!} \prod_i \frac{c_i!}{\prod_{k \neq i} c_{ik}!} \right] \cdot \left[ \prod_i \binom{m_i - e_i}{\widetilde{c_i}} \frac{\widetilde{c_i}}{\prod_{k \neq i} c_{ki}!} \right] \cdot \left[ \prod_i \prod_{k \neq i} c_{ik}! \right]
\end{equation}

\underline{Module $D$} sums over the possibilities for the number of elements in the m$'$-term of class 5. We have non-negative integers $l_{ik}$ and $r_{ik}$ with $i, k \in \{1, \cdots, r\}, i \neq k$.
\begin{equation}
    d_i = \sum_{k \neq i} l_{ik} = \sum_{k \neq i} r_{ik} \label{eq:31}
\end{equation}
\begin{equation}
    \widetilde{l_i} = \sum_{k \neq i} l_{ki} \label{eq:32}
\end{equation}
\begin{equation}
    \widetilde{r_i} = \sum_{k \neq i} r_{ki} \label{eq:33}
\end{equation}
$l_{ik}$ is the number of elements in the m$'$-term of class 5 and color $i$ and which share a row with an element of the m-term of color k. $r_{ik}$ has a similar definition with column replacing row.
\begin{equation} \label{eq:34}
    D = \sum_{r_{ik}} \sum_{l_{ik}} \prod_i \left[ \binom{m_i - e_i - \widetilde{b_i}}{\widetilde{l_i}} \cdot \frac{\widetilde{l_i}!}{\prod_{k \neq i} l_{ki}!} \right] \cdot \left[ \prod_i \binom{m_i - e_i - \widetilde{c_i}}{\widetilde{r_i}} \cdot \frac{\widetilde{r_i}!}{\prod_{k \neq i} r_{ki}!} \right] \cdot \left[ d_i! \right]
\end{equation}

We turn to the discussion of the ranges of the sums in \eqref{eq:13}. The primary variables, non-negative integers, are $m_i$, $a_i$, $e_i$, $b_{ik}$, $c_{ik}$, $l_{ik}$, and $r_{ik}$. From these are defined secondary variable $a$, $e$, $b$, $c$, $d$, $b_i$, $\widetilde{b_i}$, $c_i$, $\widetilde{c_i}$, $d_i$, $\widetilde{l_i}$, $\widetilde{r_i}$, by equations \eqref{eq:16}, \eqref{eq:17}, \eqref{eq:18}, \eqref{eq:19}, \eqref{eq:20}, \eqref{eq:25}, \eqref{eq:26}, \eqref{eq:28},
\eqref{eq:29}, \eqref{eq:31}, \eqref{eq:32}, \eqref{eq:33} respectively. Restrictions on the independence of the basic variables are given by \eqref{eq:5}, \eqref{eq:15}, and \eqref{eq:31}. In addition to these three restrictions we find the following inequalities: 
\renewcommand{\labelenumi}{(\alph{enumi})}
\begin{enumerate}
    \item From the consideration of \eqref{eq:21} we find
        \begin{equation}
            m_i + a_i + b_i + c_i + d_i \leq n \text{, each } i
        \end{equation}
    \item From the consideration of \eqref{eq:23} we find
        \begin{equation}
            a \leq n - m
        \end{equation}
    \item From the consideration of \eqref{eq:24} we find
        \begin{equation}
            e_i \leq m_{i} \text{, each } i
        \end{equation}
    \item From the consideration of \eqref{eq:27} we find
        \begin{equation}
            a + b + m \leq n
        \end{equation}
        and
        \begin{equation}
            e_i + \widetilde{b_i} \leq m_i \text{, each } i
        \end{equation}
    \item From the consideration of \eqref{eq:30} we find
        \begin{equation}
            a + c + m \leq n
        \end{equation}
        and
        \begin{equation}
            e_i + \widetilde{c_i} \leq m_i \text{, each } i
        \end{equation}
    \item From the consideration of \eqref{eq:34} we find
        \begin{equation}
            e_i + \widetilde{b_i} + \widetilde{l_i} \leq m_i \text{, each } i
        \end{equation}
        and
        \begin{equation}
            e_i + \widetilde{c_i} + \widetilde{r_i} \leq m_i \text{, each } i
        \end{equation}
\end{enumerate}

With the definition of the modules, $M$ through $T$, and these restrictions on the sum over the primary variables, equation \eqref{eq:13} provides an exact expression for the expectation of the product of permanents.

Referring back to the developments surrounding equation \eqref{eq:7} we write the exact expression for the expectation in the form
\begin{equation}
    \sum e^S
\end{equation}
We consider $S$ as a function of the primary variables. Then the ``statistical assumption'' says that if $a$, $e$, $b$, $c$, $d$, are fixed at some large values, that is on the scale of $n$, then with these restrictions, the largest value of $S$ is obtained when
\begin{equation} \label{eq:45}
    m_i = m_j \text{, all } i, j
\end{equation}
\begin{equation}
    a_i = a_j \text{, all } i, j
\end{equation}
\begin{equation}
    e_i = e_j \text{, all } i, j
\end{equation}
\begin{equation}
    b_{ij} = b_{kl} \text{, all } i, j, k, l
\end{equation}
\begin{equation}
    c_{ij} = c_{kl} \text{, all } i, j, k, l
\end{equation}
\begin{equation}
    l_{ij} = l_{kl} \text{, all } i, j, k, l
\end{equation}
\begin{equation}
    r_{ij} = r_{kl} \text{, all } i, j, k, l
\end{equation}
\begin{equation} \label{eq:52}
    l_{ij} = r_{ij} \text{, all } i, j
\end{equation}
Relations \eqref{eq:45} through \eqref{eq:52} holding on the scale of 1.

These equalities we argue hold by statistical considerations of how the colors distribute themselves once the structure of the locations is specified.

We now set
\begin{equation} \label{eq:53}
    m_i = \frac{m}{r}
\end{equation}
\begin{equation}
    a_i = \frac{a}{r}
\end{equation}
\begin{equation}
    e_i = \frac{e}{r}
\end{equation}
\begin{equation}
    b_{ij} = \frac{b}{r(r - 1)}
\end{equation}
\begin{equation}
    c_{ij} = \frac{c}{r(r - 1)}
\end{equation}
\begin{equation} \label{eq:58}
    l_{ij} = r_{ij} = \frac{d}{r(r -1)}
\end{equation}
We also set $b = c$ for simplicity as will turn out true if we let them remain independent through the following discussion. We use the approximation \eqref{eq:3} for $\ln(z!)$ and write $S$ using \eqref{eq:53} through \eqref{eq:58} with $b = c$.
\begin{equation}
    S \cong LM + LA + LE + 2 LB + LD + LT
\end{equation}
with
\begin{equation}
    LM = -r F(n) + 2 F(n) - 2F(n - m) - r F \left( \frac{m}{r} \right)
\end{equation}
from \eqref{eq:14}, and
\begin{equation}
    LA = 2 F(n - m) - 2 F(n - m - a) - r F(\frac{a}{r})
\end{equation}
from \eqref{eq:23}, and
\begin{equation}
    LE = r F(\frac{m}{r})- r F(\frac{m}{r} - \frac{e}{r}) - r F \left( \frac{e}{r} \right)
\end{equation}
from \eqref{eq:24}, and
\begin{equation}
    LB = F(n - m - a) - F(n - m - a - b) - r(r-1) F \left( \frac{b}{r(r - 1)} \right) + r F \left( \frac{m}{r} - \frac{e}{r} \right) - r F \left( \frac{m}{r} - \frac{e}{r} - \frac{b}{r} \right)
\end{equation}
from \eqref{eq:27}, and
\begin{equation}
    LD = 2 r F \left( \frac{m - e - b}{r} \right) - 2 r(r-1) F \left( \frac{d}{r(r - 1)} \right) + r F \left( \frac{d}{r} \right) - 2 r F \left( \frac{m - e - b - d}{r} \right)
\end{equation}
from \eqref{eq:34}, and
\begin{equation}
    LT = r F \left( n - \frac{m + a + 2b + d}{r} \right)
\end{equation}
from \eqref{eq:21}.

Using Lagrange multipliers to impose equation \eqref{eq:15} we write
\begin{equation}
    SS = S - l \cdot (a + e + 2b + d)
\end{equation}
and attempt to solve the five following equations together.
\begin{equation}
    \frac{\partial SS}{\partial a} = 0 \text{ or } (n - m - a - b)^2 - L \left( \frac{a}{r} \right)\cdot  \left( n - \frac{m + a + 2b + d}{r} \right) = 0
\end{equation}
\begin{equation}
    \frac{\partial SS}{\partial b} = 0 \text{ or } (n - m - a - b) \cdot\left( \frac{m - e - b - d}{r} \right) - L \left( \frac{b}{r(r - 1)} \right) \cdot \left( n - \frac{m + a + 2b + d}{r} \right) = 0
\end{equation}
\begin{equation}
    \frac{\partial SS}{\partial d} = 0 \text{ or } \left( \frac{(m - e - b - d)}{r} \right)^2 - L d \frac{1}{(r(r - 1)^2)} \cdot \left( n - \frac{m + a + 2b + d}{r} \right) = 0
\end{equation}
\begin{equation}
    \frac{\partial SS}{\partial e} = 0 \text{ or } \left( \frac{(m - e - b - d)}{r} \right)^2 - L \left( \frac{m}{r} - \frac{e}{r} \right)\cdot \left( \frac{e}{r} \right) = 0
\end{equation}
\begin{equation}
    a + e + 2b + d = m'
\end{equation}

Amazingly enough Maple found the following analytic solution to this problem. We write $m = pn$ and $m' = qn$ and then have
\begin{equation}
    a = n q (1 - p)^2 \frac{r}{(r - p)}
\end{equation}
\begin{equation}
    b = n (1 - p) (r - 1) \frac{p q}{(r - p)}
\end{equation}
\begin{equation}
    d = n (r - 1)^2 p^2 \frac{q}{(r(r - p))}
\end{equation}
\begin{equation}
    e = n \frac{p q}{r}
\end{equation}
\begin{equation}
    L = (1 - q)^2 \frac{r^2}{q (r - q)}
\end{equation}
Substituting these values into $S$ we find, again amazingly
\begin{multline}
    S = n \bigg[ \Big( -p \ln p + (2p - r) \ln r + 2 (p - 1) \ln (1 - p)+ (r - p) \ln (r - p) \Big)\\
    + \Big( -q \ln q + (2q - r) \ln r + 2 (q - 1) \ln (1 - q) + (r - q) \ln (r - q) \Big) \bigg]
\end{multline}

\end{document}